\documentclass[10pt]{amsart}

\usepackage{amssymb,amsmath,amsthm,amsfonts}
\usepackage[all]{xy}

\theoremstyle{plain}
\newtheorem{Theorem}{Theorem}
\newtheorem{Proposition}[Theorem]{Proposition}
\newtheorem{Lemma}[Theorem]{Lemma}
\newtheorem{Corollary}[Theorem]{Corollary}
\newtheorem{Question}[Theorem]{Question}
\newtheorem{Observation}[Theorem]{Observation}

\theoremstyle{definition}
\newtheorem{Remark}[Theorem]{Remark}
\newtheorem{Definition}[Theorem]{Definition}

\begin{document}
\title{On direct product subgroups of $\mathrm{SO}_3(\mathbb{R})$}
\author{Diego Rattaggi}
\thanks{Supported by the Swiss National Science Foundation, No.\ PP002--68627}
\email{rattaggi@hotmail.com}
\date{\today}
\begin{abstract}
Let $G_1 \times G_2$ be a subgroup of $\mathrm{SO}_3(\mathbb{R})$
such that the two factors $G_1$ and $G_2$ are non-trivial groups.
We show that if $G_1 \times G_2$ is not abelian,
then one factor is the (abelian) group of order $2$,
and the other factor is non-abelian and contains an element of order $2$.
There exist finite and infinite such non-abelian subgroups.
\end{abstract}
\maketitle

Let $F_2$ be the free group of rank $2$. It is well-known that 
the group $\mathrm{SO}_3(\mathbb{R})$ has subgroups isomorphic to $F_2$, e.g.\
\[
\left\langle 
\left( \begin{array}{ccc}
1  &  0            &  0 \\
0  &  -3/5 &  -4/5 \\
0  &   \phantom{-}4/5 &  -3/5 \notag
\end{array}
\right), \;
\left( \begin{array}{ccc}
-3/5 & 0 & \phantom{-}4/5 \\
0 & 1 & 0 \\
-4/5 & 0 & -3/5 \notag
\end{array}
\right)
\right\rangle_{\mathrm{SO}_3(\mathbb{R})} \cong F_2,
\]
and subgroups isomorphic to $\mathbb{Z} \times \mathbb{Z}$, like
\[
\left\langle 
\left( \begin{array}{ccc}
1  &  0            &  0 \\
0  &  -3/5 &  -4/5 \\
0  &   \phantom{-}4/5 &  -3/5 \notag
\end{array}
\right), \;
\left( \begin{array}{ccc}
1  &  0            &  0 \\
0  &  -15/17 &  -8/17 \\
0  &   \phantom{-}8/17 &  -15/17 \notag
\end{array}
\right)
\right\rangle_{\mathrm{SO}_3(\mathbb{R})} \cong \mathbb{Z} \times \mathbb{Z}.
\]
However, $\mathrm{SO}_3(\mathbb{R})$ has no subgroups isomorphic
to $\mathbb{Z} \times F_2$. More precisely,
if $G_1 \times G_2$ is a non-abelian subgroup of $\mathrm{SO}_3(\mathbb{R})$ 
such that $G_1$, $G_2$ are non-trivial,
then $G_1$, $G_2$ both contain an element of order $2$,
and moreover $G_1$ or $G_2$ is abelian.
We will give an elementary proof of these results 
(Proposition~\ref{Prop7} and Proposition~\ref{Prop14}) 
using the Hamilton quaternion algebra $\mathbb{H}(\mathbb{R})$.
Additionally, we will show in Proposition~\ref{Prop16}
that any non-trivial element in the abelian factor has order $2$
and in Theorem~\ref{Thm18} that in fact the abelian factor is the group of order $2$.

Recall that elements $x \in \mathbb{H}(\mathbb{R})$ are of the form
$x = x_0 + x_1 i + x_2 j + x_3 k$, where $x_0, x_1, x_2, x_3 \in \mathbb{R}$,
and multiplication in $\mathbb{H}(\mathbb{R})$ is induced by the rules 
$i^2 = j^2 = k^2 = -1$ and $ij = -ji = k$. The \emph{norm} of $x$ is by definition
$|x|^2 = x_0^2 + x_1^2 +x_2^2 +x_3^2 \in \mathbb{R}$. 
We say that $x, y \in \mathbb{H}(\mathbb{R})$
are \emph{perpendicular} (denoted by $x \perp y$),
if $x_1 y_1 + x_2 y_2 + x_3 y_3 = 0$
(i.e.\ if $(x_1,x_2,x_3)^T$, $(y_1,y_2,y_3)^T$ are
perpendicular as vectors in $\mathbb{R}^3$).
There is a surjective homomorphism
$\vartheta$ from the multiplicative group $\mathbb{H}(\mathbb{R}) \setminus \{ 0 \}$
to $\mathrm{SO}_3(\mathbb{R})$ defined by
\[
x \mapsto 
\frac{1}{|x|^2} \left(
\begin{array}{ccc}
x_0^2 + x_1^2 - x_2^2 - x_3^2  & 2(x_1 x_2 - x_0 x_3)           & 2(x_1 x_3 + x_0 x_2) \\
2(x_1 x_2 + x_0 x_3)           & x_0^2 - x_1^2 + x_2^2 - x_3^2  & 2(x_2 x_3 - x_0 x_1) \\
2(x_1 x_3 - x_0 x_2)           & 2(x_2 x_3 + x_0 x_1)           & x_0^2 - x_1^2 - x_2^2 + x_3^2 \notag
\end{array}
\right).
\]
It is easy to check that 
\[
\mathrm{ker}(\vartheta) = Z(\mathbb{H}(\mathbb{R}) \setminus \{ 0 \}) =
\{ x \in \mathbb{H}(\mathbb{R}) \setminus \{ 0 \} : x = x_0 \}
\]
which we will identify with $\mathbb{R} \setminus \{ 0 \}$.
Note that if $x \in \mathbb{H}(\mathbb{R}) \setminus \mathbb{R}$,
then the axis of the rotation $\vartheta(x)$ is the line through 
$(0,0,0)^T$ and $(x_1,x_2,x_3)^T$ in $\mathbb{R}^3$.
Next, we prove three basic lemmas about (anti-)commutation of
quaternions.

\begin{Lemma} \label{Lemma1}
Let $x, y \in  \mathbb{H}(\mathbb{R}) \setminus \{ 0 \}$. 
Then $xy = -yx$, if and only if $x_0 = y_0 = 0$ and $x \perp y$. 
\end{Lemma}

\begin{proof}
Only using quaternion multiplication,
we get $xy = -yx$ if and only if the following four equations hold:
\begin{align}
x_1 y_1 + x_2 y_2 + x_3 y_3 &= x_0 y_0 \notag \\ 
x_0 y_1 + x_1 y_0 &= 0 \notag \\
x_0 y_2 + x_2 y_0 &= 0 \notag \\
x_0 y_3 + x_3 y_0 &= 0. \notag
\end{align}

Thus if $x_0 = y_0 = 0$ and $x \perp y$, then clearly $xy = -yx$.

To prove the converse, suppose that $xy = -yx$
and (by contradiction) $x_0 \ne 0$.
Then from the four equations, we have
$x_0 y_0 - x_1 y_1 - x_2 y_2 - x_3 y_3 = 0$ and
\[
y_1 = \frac{-x_1 y_0}{x_0}, \quad y_2 = \frac{-x_2 y_0}{x_0}, \quad y_3 = \frac{-x_3 y_0}{x_0}.
\]
It follows that
\[
x_0 y_0 + \frac{x_1^2 y_0}{x_0} + \frac{x_2^2 y_0}{x_0}  + \frac{x_3^2 y_0}{x_0} = 0
\] 
and therefore $y_0 |x|^2 = 0$. Since $|x|^2 \geq x_0^2 > 0$, we get $y_0 = 0$
which implies $y_1 = 0$, $y_2 = 0$ and $y_3 = 0$, hence the contradiction $y = 0$,
and we conclude $x_0 =0$.
The four original equations become
$x_1 y_1 + x_2 y_2 + x_3 y_3 = 0$ (i.e.\ $x \perp y$ as required)
and $x_1 y_0 = 0$, $x_2 y_0 = 0$, $x_3 y_0 = 0$, which implies $y_0 = 0$
(using $x \ne 0$) and we are done.
\end{proof}

\begin{Lemma} \label{Lemma2}
Two quaternions $x, y \in \mathbb{H}(\mathbb{R})$
commute, if and only if $(x_1,x_2,x_3)^T$
and $(y_1,y_2,y_3)^T$ are linearly dependent over $\mathbb{R}$.
\end{Lemma}

\begin{proof}
This follows from the computation
\begin{align}
xy - yx &=
2(x_2 y_3 - x_3 y_2)i + 2(x_3 y_1 - x_1 y_3)j + 2(x_1 y_2 - x_2 y_1)k
\notag \\
&= 2
\left| 
\begin{array}{ccc}
i & x_1 & y_1 \\
j & x_2 & y_2 \\
k & x_3 & y_3 
\end{array} 
\right|. \notag
\end{align}
\end{proof}

\begin{Lemma} \label{Lemma3}
Let $x, y, z \in \mathbb{H}(\mathbb{R}) \setminus \mathbb{R}$.
If $xy = yx$ and $xz = zx$, then $yz = zy$.
In other words, the group $\mathbb{H}(\mathbb{R}) \setminus \{ 0 \}$
is commutative transitive on non-central elements.
\end{Lemma}

\begin{proof}
By assumption we have
\[
\begin{pmatrix}
x_1 \\
x_2 \\
x_3 
\end{pmatrix},
\begin{pmatrix}
y_1 \\
y_2 \\
y_3 
\end{pmatrix},
\begin{pmatrix}
z_1 \\
z_2 \\
z_3 
\end{pmatrix}
\ne
\begin{pmatrix}
0 \\
0 \\
0 
\end{pmatrix}.
\]
The statement follows now directly from Lemma~\ref{Lemma2}.
\end{proof}

To describe the structure of direct product subgroups of
$\mathrm{SO}_3(\mathbb{R})$, we give some general definitions.

\begin{Definition}
We call a direct product $G_1 \times G_2$ \emph{non-trivial},
if both $G_1$ and $G_2$ are non-trivial groups.
\end{Definition}

\begin{Definition}
We say that the group $G$ satisfies property
\begin{itemize}
\item[$(P_1)$,] if $G$ is abelian.
\item[$(P_2)$,] if $G$ is CSA, i.e.\ if all its maximal abelian subgroups
are malnormal (in other words, if for any maximal abelian subgroup 
$H < G$ and any $g \in G \setminus H$ the intersection of 
$gHg^{-1}$ with $H$ is trivial).
\item[$(P_3)$,] if $G$ is commutative transitive, i.e.\ if
$xy = yx$, $xz = zx$ always implies $yz = zy$
(provided $x,y,z \in G \setminus \{ 1 \}$).
\item[$(P_4)$,] if any non-trivial direct product subgroup
$G_1 \times G_2 < G$ is abelian
(equi\-va\-lently, if in any non-trivial direct product subgroup
$G_1 \times G_2 < G$ both factors $G_1$, $G_2$ are abelian).
\item[$(P_5)$,] if any non-trivial direct product subgroup
$G_1 \times G_2 < G$ is abelian, or exactly one factor is the abelian group of order $2$
and the other factor is a non-abelian group containing an element of order $2$.
\item[$(P_6)$,] if any non-trivial direct product subgroup
$G_1 \times G_2 < G$ is abelian, or exactly one factor is abelian such
that the non-abelian factor contains an element of order $2$
and any non-trivial element in the abelian factor has order $2$. 
\item[$(P_7)$,] if any non-trivial direct product subgroup
$G_1 \times G_2 < G$ is abelian or both factors $G_1$, $G_2$ 
contain an element of order $2$.
\item[$(P_8)$,] if any torsion-free non-trivial direct product subgroup
$G_1 \times G_2 < G$ is abelian. 
\item[$(P_9)$,] if $G$ contains no subgroup $\mathbb{Z} \times F_2$. 
\item[$(P_{10})$,] if $G$ contains no subgroup $F_2 \times F_2$.
\item[$(R_3)$,] if $G$ is commutative transitive on non-central elements, 
i.e.\ if $xy = yx$, $xz = zx$ always implies $yz = zy$
(provided $x,y,z \in G \setminus ZG$).
\item[$(R_4)$,] if any non-trivial direct product subgroup
$G_1 \times G_2 < G$ is abelian, or one factor is 
non-abelian and the other factor is contained in the center of $G$.
\item[$(R_6)$,] if in any non-trivial direct product subgroup
$G_1 \times G_2 < G$ at least one factor is abelian.
\end{itemize}
\end{Definition}

\begin{Remark}
The arrows in the following diagram stand for implications.
For example ``$(P_1) \longrightarrow (P_2)$'' means
``if a group $G$ satisfies property $(P_1)$,
then $G$ satisfies property $(P_2)$''.
These implications follow directly from the given definitions, except 
maybe $(P_2) \longrightarrow (P_3)$ which is also easy to prove, see 
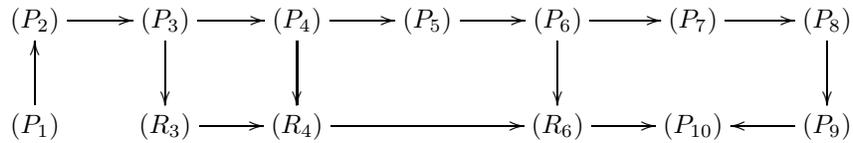
\cite[Proposition~7]{MR}.
\[
\xymatrix{
(P_2) \ar[r] & (P_3) \ar[r] \ar[d] & (P_4) \ar[r] \ar[d] & (P_5) \ar[r] & (P_6) \ar[r] \ar[d] & (P_7) \ar[r] & (P_8) \ar[d]\\
(P_1) \ar[u] & (R_3) \ar[r]        & (R_4) \ar[rr]       &              & (R_6) \ar[r]        & (P_{10})     & (P_9) \ar[l]        
}
\]
\end{Remark}

We will show in Proposition~\ref{Prop7} that
$\mathrm{SO}_3(\mathbb{R})$ satisfies property 
$(P_7)$, and in Proposition~\ref{Prop14} that
$\mathrm{SO}_3(\mathbb{R})$ satisfies property $(R_6)$,
using the map $\vartheta$ and our lemmas on quaternions.
These results will be refined in Proposition~\ref{Prop16} 
and Theorem~\ref{Thm18} to prove that $\mathrm{SO}_3(\mathbb{R})$ 
satisfies property $(P_6)$ and $(P_5)$.

For a group with trivial center, e.g.\ for $\mathrm{SO}_3(\mathbb{R})$,
properties $(P_4)$ and $(R_4)$ are equivalent.
In Observation~\ref{Obs13}, we illustrate by two examples
that $\mathrm{SO}_3(\mathbb{R})$ 
does not satisfy property $(P_4)$
(and hence does not satisfy property $(R_4)$).
As a preparation, Observation~\ref{Obs11} shows that
$\mathrm{SO}_3(\mathbb{R})$ does not satisfy property $(P_3)$.

\begin{Proposition} \label{Prop7}
The group $\mathrm{SO}_3(\mathbb{R})$ satisfies property $(P_7)$.
\end{Proposition}

\begin{proof}
Let $G_1 \times G_2$ be a non-trivial direct product subgroup
of $\mathrm{SO}_3(\mathbb{R})$ and suppose that $G_1$
or $G_2$ does not contain an element of order $2$.
We have to prove that $G_1 \times G_2$ is abelian.
Let $E$ be the identity matrix in $\mathrm{SO}_3(\mathbb{R})$, and
take any $A \in G_1 \setminus \{E\}$, $B,C \in G_2 \setminus \{E\}$. 
Then $AB = BA$ and $AC = CA$.
Take any $x, y, z \in \mathbb{H}(\mathbb{R}) \setminus \mathbb{R}$
such that $\vartheta(x) = A$, $\vartheta(y) = B$ and $\vartheta(z) = C$.
We have $\vartheta(x) \vartheta(y) = \vartheta(y) \vartheta(x)$,
hence $xyx^{-1}y^{-1} \in \mathrm{ker}(\vartheta)$, i.e.\
$xy = \lambda yx$ for some $\lambda \in \mathbb{R} \setminus \{ 0 \}$.
Taking the norm, and using the rule $|xy|^2 = |x|^2 |y|^2$,
we see that $\lambda \in \{-1, 1\}$, in other words
$xy = yx$ or $xy = -yx$.
Similarly, $AC = CA$ implies that $xz = zx$ or $xz = -zx$.

In the case $xy = -yx$, we get $x_0 = y_0 = 0$ by Lemma~\ref{Lemma1}.
But then $x^2, y^2 \in \mathbb{R} \setminus \{ 0 \}$ and
$A^2 = \vartheta(x^2) = E$, $B^2 = \vartheta(y^2) = E$,
hence both $G_1$ and $G_2$ contain an element of order $2$,
a contradiction to our assumption.
In the same way, if $xz = -zx$, then we get the contradiction $A^2 = C^2 = E$.

Hence we always have $xy = yx$ and $xz = zx$.
Using Lemma~\ref{Lemma3}, we get $yz = zy$ and therefore $BC = CB$.
This shows that $G_2$ is abelian. 
Similarly, taking two matrices in $G_1 \setminus \{E\}$
and one matrix in $G_2 \setminus \{E\}$, one shows that
$G_1$ is abelian.
\end{proof}

\begin{Corollary}
The group $\mathrm{SO}_3(\mathbb{R})$ contains no subgroup 
$\mathbb{Z} \times F_2$ and no subgroup $F_2 \times F_2$.
\end{Corollary}

\begin{proof}
Property $(P_7)$ implies property $(P_9)$ and $(P_{10})$.
\end{proof}

\begin{Remark}
A group is called \emph{coherent} if every finitely generated subgroup
is finitely presented. Any group containing a subgroup $F_2 \times F_2$
is incoherent. Therefore the non-existence of subgroups $F_2 \times F_2$
is a necessary condition for coherence, 
although it is not a sufficient condition
since there are for example incoherent (hyperbolic) groups 
(using \cite{Rips}) not containing $\mathbb{Z} \times F_2$ subgroups.
It is a question of Serre (\cite[p.734]{Serre}) whether $\mathrm{GL}_3(\mathbb{Q})$
is coherent.
\end{Remark}

\begin{Question}
Is $\mathrm{SO}_3(\mathbb{R})$ coherent?
\end{Question}

Using the idea of the proof of Proposition~\ref{Prop7},
we see that any subgroup of $\mathrm{SO}_3(\mathbb{R})$ 
which does not contain elements of order $2$
(in particular any torsion-free subgroup of $\mathrm{SO}_3(\mathbb{R})$)
is commutative transitive. However $\mathrm{SO}_3(\mathbb{R})$
itself is not commutative transitive:

\begin{Observation} \label{Obs11}
The group $\mathrm{SO}_3(\mathbb{R})$ does not satisfy property $(P_3)$.
\end{Observation}
This observation will directly follow from Observation~\ref{Obs13},
but we give a short alternative proof here.  
\begin{proof}
Take
\[
A := \left( \begin{array}{rrr}
1  &   0 &  0 \\
0  &  -1 &  0 \\
0  &   0 & -1 \notag
\end{array}
\right), \;
B := \left( \begin{array}{rrr}
1  &   0 &  0 \\
0  &   0 &  -1 \\
0  &   1 &  0 \notag
\end{array}
\right), \;
C := \left( \begin{array}{rrr}
-1  &   0 &  0 \\
0  &   1 &   0 \\
0  &   0 &  -1 \notag
\end{array}
\right),
\]
then $AB = BA$ and $AC = CA$, but $BC \ne CB$.

Note that $A = \vartheta(i)$, $B = \vartheta(1+i)$, $C = \vartheta(j)$
and $i (i+1) = (i+1) i$, $ij = -ji$, $(i+1)j \ne \pm j(i+1)$.
\end{proof}

\begin{Corollary}
There is a group $G$ which is commutative transitive on non-central
elements, but such that $G/Z(G)$ is not commutative transitive on
non-central elements (and therefore such that $G/Z(G)$ is not
commutative transitive).  
\end{Corollary}

\begin{proof}
Take $G = \mathbb{H}(\mathbb{R}) \setminus \{ 0 \}$
such that $G/ZG \cong \mathrm{SO}_3(\mathbb{R})$
and note that $Z(\mathrm{SO}_3(\mathbb{R}))$ is the trivial group.
\end{proof}

The matrices $A$, $B$, $C$ from the proof of
Observation~\ref{Obs11} generate a non-abelian subgroup 
$\langle A,B,C \rangle$ of $\mathrm{SO}_3(\mathbb{R})$.
However, this group cannot be used
to prove that $\mathrm{SO}_3(\mathbb{R})$ does not satisfy property $(P_4)$,
since $A = B^2$ and $\langle A,B,C \rangle = \langle B,C \rangle$ is the dihedral group of order $8$
which is not decomposable as a non-trivial direct product.
Nevertheless, there \emph{are} non-abelian non-trivial 
direct product subgroups of $\mathrm{SO}_3(\mathbb{R})$.

\begin{Observation} \label{Obs13}
The group $\mathrm{SO}_3(\mathbb{R})$ does not satisfy property $(P_4)$.
\end{Observation}

\begin{proof}
We give two examples of a non-abelian non-trivial 
direct product subgroup of $\mathrm{SO}_3(\mathbb{R})$, at first an
infinite example.

Let $A = \vartheta(i)$, $C = \vartheta(j)$ as in the proof of
Observation~\ref{Obs11} and let
\[
\tilde{B} := \vartheta(1+2i) = 
\left( \begin{array}{ccc}
1  &  0            &  0 \\
0  &  -3/5 &  -4/5 \\
0  &   \phantom{-}4/5 &  -3/5 \notag
\end{array}
\right).
\]
We claim that $\langle A,\tilde{B},C \rangle$ is a non-abelian non-trivial 
direct product subgroup of $\mathrm{SO}_3(\mathbb{R})$.

First we want to show by contradiction that 
$A \notin \langle \tilde{B},C \rangle$.
Since $C \tilde{B} = \tilde{B}^{-1} C$ and 
$C \tilde{B}^{-1} = \tilde{B} C$,
any word in the letters $\tilde{B}$, $\tilde{B}^{-1}$, $C = C^{-1}$
can be brought to the form $\tilde{B}^n C$ or $\tilde{B}^n$ for some
$n \in \mathbb{Z}$.
If we suppose that $A \in \langle \tilde{B},C \rangle$,
then, looking at the upper left entry 
(which is $1$ in $A$ and $\tilde{B}$, but $-1$ in $C$), 
we see that $A$ cannot be written
as $\tilde{B}^n C$ and therefore $A = \tilde{B}^n$ for some $n \in \mathbb{Z} \setminus \{ 0 \}$.
But since $A$ has order $2$, we get $\tilde{B}^{2n} = E$, which
contradicts the fact that $\tilde{B}$ has infinite order.

Since $\langle A \rangle$ has only two elements and $A \notin \langle \tilde{B},C \rangle$,
we get $\langle A \rangle \cap \langle \tilde{B}, C \rangle = \{ E \}$.
Moreover, it is easy to check that $A$ commutes with $\tilde{B}$ and with $C$. 
Therefore $\langle A, \tilde{B}, C \rangle < \mathrm{SO}_3(\mathbb{R})$
is a direct product of the group $\langle A \rangle$ of order $2$
and the (infinite) non-abelian (solvable) group $\langle \tilde{B}, C \rangle$.

As a finite example we can take the dihedral group of order $12$,
generated for example by the two matrices
\[
\left( \begin{array}{ccc}
1  &  0            &  0 \\
0  &  1/2 &  -\sqrt{3}/2 \\
0  &  \sqrt{3}/2 &  1/2 \notag
\end{array}
\right) \text{ and }
\left( \begin{array}{rrr}
-1  &   0 &  0 \\
0  &   1 &   0 \\
0  &   0 &  -1 \notag
\end{array}
\right).
\]
This group is isomorphic to a direct product of the (non-abelian) dihedral group of order~$6$
(which is isomorphic to the symmetric group $S_3$) 
and the group of order~$2$.
\end{proof}

\begin{Proposition} \label{Prop14}
The group $\mathrm{SO}_3(\mathbb{R})$ satisfies property $(R_6)$.
\end{Proposition}

\begin{proof}
Suppose by contradiction that $G_1 \times G_2$ is a non-trivial direct 
subgroup of $\mathrm{SO}_3(\mathbb{R})$ such that 
$G_1$ and $G_2$ are non-abelian.
First take $A,B \in G_1 \setminus \{E \}$ such that $AB \ne BA$
and $C,D \in G_2 \setminus \{E \}$ such that $CD \ne DC$.
Now take $x,y,z,w \in \mathbb{H}(\mathbb{R}) \setminus \mathbb{R}$
such that $\vartheta(x) = A$, $\vartheta(y) = B$, $\vartheta(z) = C$,
$\vartheta(w) = D$.
Then we have $xy \ne \pm yx$, $zw \ne \pm wz$
and (by the same argument as in the proof of Proposition~\ref{Prop7})
$xz = \pm zx$, $xw = \pm wx$, $yz = \pm zy$, $yw = \pm wy$.

Suppose that $xz = zx$.
If $xw = wx$ then we get by Lemma~\ref{Lemma3} the contradiction $zw = wz$,
hence $xw = -wx$. But then by Lemma~\ref{Lemma1}, 
$w_0 = 0$ and $x \perp w$. 
Since $(x_1,x_2,x_3)^T$ and $(z_1,z_2,z_3)^T$ are 
linearly dependent by Lemma~\ref{Lemma2}, we conclude $z \perp w$.
Since $wz \ne -zw$, we have $z_0 \ne 0$ by Lemma~\ref{Lemma1},
hence $yz = zy$ again by Lemma~\ref{Lemma1},
and $xy = yx$ by Lemma~\ref{Lemma3}, a contradiction.

We have shown that $xz = -zx$. Similarly, it follows that
$xw = -wx$, $yz = -zy$ and $yw = -wy$.
Lemma~\ref{Lemma1} implies $x \perp z$ and $x \perp w$.
Since $zw \ne wz$, $z$ and $w$ are linearly independent 
by Lemma~\ref{Lemma2} and span the plane perpendicular to $x$.
We also have $y \perp z$ and $y \perp w$ by Lemma~\ref{Lemma1},
hence $x$ and $y$ are linearly dependent and we get
the contradiction $xy = yx$ by Lemma~\ref{Lemma2}.
\end{proof}

\begin{Lemma} \label{Lemma15}
Let $A \in \mathrm{SO}_3(\mathbb{R})$ be a rotation of order at least $3$.
Then the centralizer of $A$ in $\mathrm{SO}_3(\mathbb{R})$ consists of all rotations about the axis
of $A$.
\end{Lemma}

\begin{proof}
Without loss of generality, we may assume that $A$ is a rotation 
of order at least $3$ about the $x$-axis, hence
\[
A = \left( \begin{array}{ccc}
1  &  0            &  0 \\
0  &  \cos \phi &  -\sin \phi \\
0  &  \sin \phi &  \phantom{-}\cos \phi \notag
\end{array}
\right),
\]
such that $\sin \phi \ne 0$.
Suppose that the matrix
\[
B = 
\left( \begin{array}{ccc}
b_1 & b_2 & b_3 \\
b_4 & b_5 & b_6 \\
b_7 & b_8 & b_9 \notag
\end{array}
\right) \in \mathrm{SO}_3(\mathbb{R})
\]
commutes with $A$.
Then $AB = BA$ gives the conditions
\begin{align}
b_9 \sin \phi &= b_5 \sin \phi \notag \\
-b_8 \sin \phi &= b_6 \sin \phi \notag
\end{align}
and
\begin{align}
b_2 (1 - \cos \phi) &= b_3 \sin \phi \notag \\
-b_3 (1 - \cos \phi) &= b_2 \sin \phi \notag \\
-b_4 (1 - \cos \phi) &= b_7 \sin \phi \notag \\
b_7 (1 - \cos \phi) &= b_4 \sin \phi. \notag
\end{align}
The first two equations imply $b_5 = b_9$ and $b_6 = -b_8$.
The third and fourth equation imply
\[
b_2 = \frac{-b_3 (1 - \cos \phi)}{\sin \phi} \, \text{ and } \,
\frac{-b_3 (1 - \cos \phi)^2}{\sin \phi} = b_3 \sin \phi,
\]
hence
\[
-b_3(1 - 2\cos \phi) = b_3(\sin^2 \phi + \cos^2 \phi) = b_3.
\]
If $b_3 \ne 0$ then $1 - 2\cos \phi = -1$, hence $\cos \phi = 1$
and we get the contradiction $\sin \phi = 0$. Thus $b_3 = 0$ and $b_2 = 0$.
Similarly, the fifth and sixth equation lead to
$b_4 = b_7 = 0$, hence
\[
B = 
\left( \begin{array}{ccc}
b_1 & 0 & 0 \\
0 & b_5 & -b_8 \\
0 & b_8 & \phantom{-}b_5 \notag
\end{array}
\right) 
\]
We exclude the case $b_1 = -1$ computing the determinant of $B$,
and conclude
\[
B
= \left( \begin{array}{ccc}
1 & 0 & 0 \\
0 & \cos \psi & -\sin \psi \\
0 & \sin \psi & \phantom{-}\cos \psi \notag
\end{array}
\right)
\]
for some $\psi$.
\end{proof}

\begin{Proposition} \label{Prop16}
The group $\mathrm{SO}_3(\mathbb{R})$ satisfies property $(P_6)$.
\end{Proposition}

\begin{proof}
Let $G_1 \times G_2$ be a subgroup of $\mathrm{SO}_3(\mathbb{R})$ 
such that $G_2$ is non-abelian and $G_1$ is abelian and non-trivial.
Using Proposition~\ref{Prop7} and Proposition~\ref{Prop14},
it remains to prove that any non-trivial element of $G_1$ has order $2$.
Therefore suppose that $A \in G_1 \setminus \{ E \}$
has order at least $3$. Then by Lemma~\ref{Lemma15}, any
element in $G_2$ is a rotation about the axis of $A$,
which contradicts our assumption that $G_2$ is non-abelian.
\end{proof}

\begin{Lemma} \label{Lemma17}
The two matrices
\[
\left( \begin{array}{ccc}
-1  &  0            &  0 \\
0  &  \cos \phi &  \phantom{-}\sin \phi \\
0  &  \sin \phi &  -\cos \phi \notag
\end{array}
\right), 
\left( \begin{array}{ccc}
-1  &  0            &  0 \\
0  &  \cos \psi &  \phantom{-}\sin \psi \\
0  &  \sin \psi &  -\cos \psi \notag
\end{array}
\right) \in \mathrm{SO}_3(\mathbb{R})
\]
commute, if and only if 
\[
\frac{\phi}{2} - \frac{\psi}{2} \in \{k \cdot \frac{\pi}{2} : k \in \mathbb{Z}\}.
\]
In particular, these two $180^{\circ}$-rotations commute, if and only if their
axes (which lie in the $yz$-plane) are identical or perpendicular.
\end{Lemma}

\begin{proof}
Matrix multiplication gives the condition 
$\sin \phi \cdot \cos \psi = \cos \phi \cdot \sin \psi$,
hence
\[
0 = \sin \phi \cdot \cos \psi - \cos \phi \cdot \sin \psi = \sin (\phi - \psi)
\] 
and
\[
\phi - \psi \in \{k \cdot \pi : k \in \mathbb{Z}\}.
\]
\end{proof}

\begin{Theorem} \label{Thm18}
The group $\mathrm{SO}_3(\mathbb{R})$ satisfies property $(P_5)$.
\end{Theorem}

\begin{proof}
Let $G_1 \times G_2$ be a subgroup of $\mathrm{SO}_3(\mathbb{R})$ 
such that $G_2$ is non-abelian and $G_1$ is abelian and non-trivial.
Applying Proposition~\ref{Prop16}, it remains to show
that $G_1$ has order $2$. 
Let $A \in G_1 \setminus \{E\}$.
Without loss of generality we may assume that $A$ is a rotation about 
the $x$-axis. It has order $2$ by Proposition~\ref{Prop16}, 
hence 
\[
A = 
\left( \begin{array}{rrr}
1  &   0 &  0 \\
0  &  -1 &  0 \\
0  &   0 & -1 \notag
\end{array}
\right).
\]
Any element in $G_1 \setminus \{E\}$ has order $2$ and commutes with $A$.
An easy computation shows that if an element in $\mathrm{SO}_3(\mathbb{R})$
commutes with $A$, then it has either the form
\[
\left( \begin{array}{ccc}
1  &  0            &  0 \\
0  &  \cos \varphi &  -\sin \varphi \\
0  &  \sin \varphi &  \phantom{-}\cos \varphi \notag
\end{array}
\right) \, \text{ or } \, 
\left( \begin{array}{ccc}
-1  &  0            &  0 \\
0  &  \cos \phi &  \phantom{-}\sin \phi \\
0  &  \sin \phi &  -\cos \phi \notag
\end{array}
\right),
\]
i.e.\ it is either a rotation about the $x$-axis, or 
a rotation about an axis in the $yz$-plane 
by an angle of $180^{\circ}$.
The only element of order $2$ of the first form is $A$ itself.
Hence if $G_1 \setminus \{ E, A\}$ is not empty, then it contains only
elements of the second form. 
Since $G_1$ is abelian, $G_1 \setminus \{ E, A\}$ 
contains by Lemma~\ref{Lemma17} at most two elements,
and $G_1$ has therefore at most $4$ elements.
However, we know by Proposition~\ref{Prop7} that also $G_2$ contains an element of order $2$
commuting with $A$, hence $G_1$ has less than $4$ elements.
Since $A \in G_1$ has order $2$, we conclude that 
$G_1$ has exactly $2$ elements.
\end{proof}

\begin{Remark}
All statements in this article remain true if we replace $\mathbb{R}$ by $\mathbb{Q}$.
The only construction where we have used irrational numbers
was in the second part of Observation~\ref{Obs13}.
\end{Remark}

\end{document}